\documentclass[11pt,a4paper]{article} 

\usepackage{amssymb,amsfonts,amsthm, amscd}
\usepackage{amsmath}
\usepackage{hyperref}

\theoremstyle{plain}
\numberwithin{equation}{section}
\newtheorem{theorem}{Theorem}[section]

\theoremstyle{remark}

\theoremstyle{definition}

\title{A new kind of uniqueness theorems for inverse Sturm-Liouville problems}
\author{Yuri Ashrafyan}

\begin{document}

\maketitle

\begin{abstract}
We prove Marchenko-type uniqueness theorems for inverse Sturm-Liouville problems. Moreover, we prove a generalization of Ambarzumyan theorem.
\end{abstract}


\section{Introduction}\label{sec1}

Let us denote by $L(q, \alpha, \beta)$ the Sturm-Liouville boundary value problem
\begin{gather}
-y''+q(x)y=\mu y,\quad x\in (0, \pi),\ \mu\in \mathbb{C},\label{eq1.1}\\
y(0)\cot \alpha + y'(0) = 0,\qquad \alpha\in (0, \pi),\label{eq1.2}\\
y(\pi)\cot \beta + y'(\pi)  = 0,\qquad \beta\in (0, \pi),\label{eq1.3}
\end{gather}
{where $q$ is a real-valued, summable function, $q \in L^1_\mathbb{R}(0, \pi)$.
At the same time, $L(q, \alpha, \beta)$ denotes the self-adjoint operator, generated by
problem \eqref{eq1.1}-\eqref{eq1.3} (see, e.g. \cite{Naimark, Marchenko1, LevSar}).
It is known, that under the above conditions the spectrum of operator
 $L(q, \alpha, \beta)$ is discrete and consists of real, simple eigenvalues (see, e.g. \cite{Marchenko1}, \cite{Yurko1}), which we denote by $\mu_n=\mu_n(q,\alpha,\beta)$, $n \geq 0$, emphasizing the dependence of $\mu_n$ on $q$, $\alpha$
and $\beta$. We assume that eigenvalues are enumerated in the increasing order, i.e.,
$$
\mu_0(q, \alpha, \beta) < \mu_1(q, \alpha, \beta) < \dots < \mu_n(q, \alpha, \beta) < \dots.
$$

Let $\varphi(x,\mu)$ be a solution of equation \eqref{eq1.1}, which satisfies the initial conditions
\begin{equation}\label{eq1.4}
\varphi(0,\mu)=1,\quad \varphi'(0,\mu)=-\cot\alpha.
\end{equation}
The eigenvalues $\mu_n=\mu_n(q, \alpha, \beta)$, $n \geq 0$, of
$L(q, \alpha, \beta)$ are the solutions of equation
\begin{equation*}
\varphi(\pi, \mu)\cot \beta+\varphi'(\pi, \mu)=0. 
\end{equation*}
It is easy to see that the functions $\varphi(x, \mu_n)$, $n \geq 0$, are the eigenfunctions, corresponding to the eigenvalue $\mu_n$.
The squares of the $L^2$-norm of these eigenfunctions:
\begin{equation*}
a_n = a_n(q,\alpha,\beta) :=\int_0^{\pi} |\varphi(x, \mu_n)|^2 dx,\quad    n \geq 0,
\end{equation*}
are called norming constants.
The eigenvalues and norming constants are called spectral data (besides these, there are other quantities, which are also called spectral data).
The inverse Sturm-Liouville problem is to reconstruct the quantities  $q, \ \alpha, \ \beta$ by some spectral data.

Let $L=L(q, \alpha, \beta)$ and $L_0=L(q_0, \alpha_0, \beta_0)$ be two operators.
The following assertion is usually called uniqueness theorem of Marchenko\footnote{ The theorem of Marchenko is more general, see e.g. \cite{Marchenko2}--\cite{FreYur}.}.
\begin{theorem}[Marchenko \cite{Marchenko2}]\label{thm1.1}
Let  $q \in L^1_\mathbb{R}(0, \pi)$.
If
\begin{gather}
\mu_n(q, \alpha, \beta) = \mu_n(q_0, \alpha_0, \beta_0), \label{eq1.5}\\
a_n(q, \alpha, \beta) = a_n(q_0, \alpha_0, \beta_0), \label{eq1.6}
\end{gather}
for all $n \geq 0$, then $\alpha = \alpha_0, \ \beta = \beta_0$ and $q(x) = q_0(x)$ almost everywhere.
\end{theorem}

One of the results of the present paper is the following theorem, which, in some sense, is a generalization of Marchenko's uniqueness theorem.
\begin{theorem}\label{thm1.2}
	Let $q' \in L^2_\mathbb{R}(0, \pi)$.
	If
	\begin{gather}
		\mu_n(q, \alpha_0, \beta) = \mu_n(q_0, \alpha_0, \beta_0), \label{eq1.7} \\
		a_n(q, \alpha_0, \beta) \geq a_n(q_0, \alpha_0, \beta_0), \label{eq1.8}
	\end{gather}
	for all $n \geq 0$, then $\beta = \beta_0$ and $q(x) \equiv q_0(x)$.
\end{theorem}

This kind of uniqueness theorem has not been considered before. 
The main difference between Theorems \ref{thm1.1} and \ref{thm1.2} is that the equality in \eqref{eq1.6} we replace with inequality in \eqref{eq1.8}.
Note, we assume  $q' \in L^2_\mathbb{R}(0, \pi)$ instead of general $q \in L^1_\mathbb{R}(0, \pi)$, since our proof is based on the results of Jodeit and Levitan (see \cite{JodLev}).
And the parameter $\alpha$ of boundary condition is in advance fixed $\alpha = \alpha_0$.

\textbf{Remark 1.} Some analogues of Theorem \ref{thm1.2} will be stated in Section \ref{sec5}.

\

Historically, the first work in the theory of inverse spectral problems for Sturm-Liouville operators belongs to Ambarzumyan \cite{Ambarzumyan}. 
He proved that if the eigenvalues of Sturm-Liouville operator with Neumann boundary conditions are $n^2$, then the potential $q$ is $0$ on $[0, \pi]$.
It's known that the eigenvalues $\mu_n(0, \pi/2, \pi/2)$ of operator $L(0, \pi/2, \pi/2 )$ are $n^2, \ n \geq 0$. 
The classical Ambarzumyan's theorem in our notations will be as follows.
\begin{theorem}[Ambarzumyan \cite{Ambarzumyan}]\label{thm1.3} \ 
	
	If $\mu_n(q, \pi/2, \pi/2) = \mu_n(0, \pi/2, \pi/2)=n^2$, for all $n \geq 0$, then $q(x) \equiv 0$.
\end{theorem}

This is an exception, as in general additional information is needed in order to reconstruct the potential $q$ uniquely.
There are many generalizations of Ambarzumyan's theorem in various directions, we mention several of them (see, e.g. \cite{Kuznezov, ChakAch, ChernLawWang1, ChernLawWang2, YangHuangYang, YangWang, Yurko2, YilKoy} and references therein).

Our generalization of Ambarzumyan's theorem is as follows.
\begin{theorem}\label{thm1.4} 
	Let $q' \in L^2_\mathbb{R}(0, \pi)$.
	
	If $\mu_n(q, \alpha, \pi - \alpha) = \mu_n(0, \alpha, \pi - \alpha)$, for all $n \geq 0$,
	then $q(x) \equiv 0$.
\end{theorem}

We think that Theorem \ref{thm1.4} is a natural generalization, because we use only one spectrum to reconstruct the potential $q$, without any additional conditions, as it is in the classical result.

\

\section{Preliminaries}\label{sec2}

Two operators $L=L(q, \alpha, \beta)$ and $L_0=L(q_0, \alpha_0, \beta_0)$ are called isospectral, if they have the same spectra, i.e. $\mu_n (q, \alpha, \beta) = \mu_n (q_0, \alpha_0, \beta_0), \ n \geq 0$.	
In what follows, if a certain symbol $\gamma$ denotes an object related to $L$, then $\gamma_0$ (or $\gamma^0$, depending on situation) will denote a similar object related to $L_0$.

The problem of describing all the operators $L$ isospectral with $L_0$ first was considered by Trubowitz et al. (see \cite{IsaTru, IsaMcTru, DahTru,  PosTru}) for $q \in L^2_\mathbb{R}(0, \pi)$.
The same problem was considered by Jodeit and Levitan in \cite{JodLev} for $q$, such that $q' \in L^2_\mathbb{R}(0, \pi)$.
For this aim the Gelfand-Levitan integral equation and transformation operators were used in \cite{JodLev}.
They construct the kernel $F(x,y)$ of the integral equation as follows.
Let $c_n, \ n \geq 0$, be arbitrary real numbers, converging to zero, as $n \rightarrow \infty$, so rapidly, that the function

\begin{equation}\label{eq2.1}
  F(x, y) = \sum_{n=0}^{\infty} c_n \varphi_0(x,\mu_n^0) \varphi_0(y, \mu_n^0)
\end{equation}
is continuous and all the second order partial derivatives are also continuous.
The integral equation

\begin{equation}\label{eq2.2}
  K(x, y) + F(x, y) + \int^x_0 K(x, t) F(t, y) dt = 0,\qquad 0 \leq y \leq x \leq \pi,
\end{equation}
is called Gelfand-Levitan integral equation\footnote{
Here $F(x, y)$ is a kernel of integral equation \eqref{eq2.2}, where $x$ is a parameter, $F(x, y)$ is known function and $K(x, y)$ is unknown function, as functions of $y$.}.

They proved, that if $1 + c_n a_n^0 > 0$, for all $n \geq 0$, then the integral equation \eqref{eq2.2} has a unique solution $K(x,y)$ and the function
\begin{equation*}
  \varphi(x, \mu) = \varphi_0(x, \mu) + \int_{0}^{x} K(x, t) \varphi_0(t, \mu) dt
\end{equation*}
is a solution of the differential equation \eqref{eq1.1}, with potential function

\begin{equation}\label{eq2.3}
  q(x) = q_0(x) + 2 \cfrac{d}{dx} K(x, x),
\end{equation}
and $\varphi(x,\mu)$ satisfies the initial conditions 
\begin{equation*}
\varphi(0, \mu) = 1, \qquad \varphi'(0, \mu) = -\cot \alpha,
\end{equation*}
where
\begin{equation}\label{eq2.4}
\cot \alpha = \cot \alpha_0 + \sum_{n=0}^{\infty} c_n.
\end{equation} 
It means, that the function $\varphi(x, \mu)$ satisfies the boundary condition \eqref{eq1.2} for all $\mu \in \mathbb{C}$.

Find $\beta \in (0,\pi)$, such that $\mu_n(q, \alpha, \beta) = \mu_n(q_0, \alpha_0, \beta_0)$, for all $n \geq 0$, i.e. $\varphi(x, \mu)$ should satisfy, at the point $x = \pi$, the boundary condition \eqref{eq1.3}
\[
\varphi(\pi, \mu_n^0) \cot \beta + \varphi'(\pi, \mu_n^0)=0,
\]
for this $\beta \in (0, \pi)$.
Such $\beta$ (in \cite{JodLev}) is being defined from the following relation
\begin{equation}\label{eq2.5}
  \cot \beta = \cot \beta_0 + \sum_{n=0}^{\infty} \cfrac{c_n \varphi_0^2(\pi, \mu_n^0)}{1 + c_n a_n^0}.
\end{equation}

Thus Jodeit and Levitan showed, that each admissible sequence $\{c_n\}_{n=0}^\infty$ generate an isospectral operator $L(q, \alpha, \beta)$, where $q, \ \alpha$ and $\beta$ are given by the formulae \eqref{eq2.3}, \eqref{eq2.4} and \eqref{eq2.5} respectively.
In this way they obtained all the potentials $q$, with $q' \in L^2(0, \pi)$, having a given spectrum $\mu_n^0 = \mu_n(q_0, \alpha_0, \beta_0), \ n \geq 0$.

\

\section{Proof of Theorem \ref{thm1.2}}\label{sec3}

 Consider operators $L_0=L(q_0, \alpha_0, \beta_0)$ and $L=L(q, \alpha_0, \beta)$, with the set of norming constants $a_n^0 = a_n(q_0, \alpha_0, \beta_0)$ and $a_n = a_n(q, \alpha_0, \beta)$, $n \geq 0$, respectively.
 It is known (see, e.g. \cite{JodLev}), that in this case the kernel $F(x, y)$ of the integral equation \eqref{eq2.2} is
 \begin{equation}\label{eq3.1}
     F(x, y) = \sum_{n=0}^{\infty} \left( \cfrac{1}{a_n} - \cfrac{1}{a_n^0} \right) \varphi_0(x,\mu_n^0) \varphi_0(y, \mu_n^0).
 \end{equation}
Since by the condition of Theorem \ref{thm1.2} the operators $L$ and $L_0$ are isospectral, then the formulae \eqref{eq2.3}--\eqref{eq2.5} are held.
If we compare the kernels \eqref{eq2.1} and \eqref{eq3.1}, we'll refer, that $c_n = \cfrac{1}{a_n} - \cfrac{1}{a_n^0}$.
So the formulae \eqref{eq2.4} and \eqref{eq2.5} will become
\begin{equation}\label{eq3.2}
  \cot \alpha = \cot \alpha_0 + \sum_{n=0}^{\infty} \left( \cfrac{1}{a_n} - \cfrac{1}{a_n^0} \right),
\end{equation}
\begin{equation}\label{eq3.3}
  \cot \beta = \cot \beta_0 + \sum_{n=0}^{\infty} ( a_n^0 - a_n ) \cfrac{ \varphi^2_0(\pi, \mu_n^0)}{(a_n^0)^2}.
\end{equation}
Thus, we have all the operators $L(q, \alpha, \beta)$ isospectral with $L(q_0, \alpha_0, \beta_0)$.

We supposed, that $\alpha = \alpha_0$, then by formula \eqref{eq3.2} we have
  \begin{equation}\label{eq3.4}
      \sum_{n=0}^{\infty} \left( \cfrac{1}{a_n} - \cfrac{1}{a_n^0} \right) = 0.
  \end{equation}
Since $a_n \geq a_n^0$, for all $n \geq 0 $, thus from the equation \eqref{eq3.4} it refers that $a_n = a_n^0$, for all $n \geq 0 $.
Thus, from Marchenko uniqueness theorem\ref{thm1.1} we obtain  $q(x) \equiv q_0(x)$ and $\beta = \beta_0$.

This completes the proof.

\textbf{Remark 2.} From the equation \eqref{eq3.4} it follows, that the condition $a_n \geq a_n^0$ can be changed with $a_n \leq a_n^0$.
From the relation \eqref{eq3.3} it follows, that we can assume $\beta = \beta_0$, instead of $\alpha = \alpha_0$, with the condition $a_n \geq a_n^0$ or $a_n \leq a_n^0$ and then we will also obtain $q(x) \equiv q_0(x)$ and $\alpha = \alpha_0$. 

\

\section{Proof of Theorem \ref{thm1.4}}\label{sec4}

Consider an operator $L(q, \alpha, \pi - \alpha)$ and an even operator\footnote{A problem $L(q, \alpha, \beta)$ is said to be even, if $q(x) = q(\pi - x)$ and $\alpha + \beta = \pi$.}
$L(0, \alpha, \pi - \alpha)$.

N. Levinson proved \cite{Levinson} (see also \cite{Harutyunyan}), that an operator $L$ is even if and only if
\begin{equation}\label{eq4.1}
\varphi(\pi, \mu_n)=(-1)^n, \qquad n \geq 0.
\end{equation}

The condition of the theorem means, that the operator $L(q, \alpha, \pi - \alpha)$ is isospectral with $L(0, \alpha, \pi - \alpha)$. 
Since the method of Jodeit and Levitan has described all the isospectral operators for potential function $q$, with $q' \in L^2(0, \pi)$, then there exists sequence $\{c_n\}_{n=0}^{\infty}$, such that $1 + c_n a_n^0 > 0$, for all $n \geq 0$, and $\{c_n\}_{n=0}^{\infty}$ has the properties described in Section \ref{sec2} and the formulae \eqref{eq2.3}--\eqref{eq2.5} are held for operators $L(q, \alpha, \pi - \alpha)$ and  $L(0, \alpha, \pi - \alpha)$.

Therefore, taking into account, that $q_0(x) \equiv 0$, $\ \alpha_0 = \alpha$, $\ \beta_0 = \beta = \pi - \alpha$ and \eqref{eq4.1}, then the relations \eqref{eq2.3}--\eqref{eq2.5}, which connect these two operators, will become
\begin{gather}
q(x) = 2 \cfrac{d}{dx} K(x, x), \label{eq4.2} \\
\sum_{n=0}^{\infty} c_n = 0. \label{eq4.3} \\
\sum_{n=0}^{\infty} \cfrac{c_n}{1 + c_n a_n^0} = 0. \label{eq4.4}
\end{gather}

If we subtract \eqref{eq4.3} from \eqref{eq4.4} we will obtain
\begin{equation}\label{eq4.8}
	\sum_{n=0}^{\infty} \cfrac{c_n^2 a_n^0}{1 + c_n a_n^0} = 0.
\end{equation}
Since $1 + c_n a_n^0 > 0$ and $a_n^0 > 0$, for all $n \geq 0$, then from the equation \eqref{eq4.8} we obtain, that $c_n = 0, \ n \geq 0$.
Thus, from the equations \eqref{eq2.1}, \eqref{eq2.2} and \eqref{eq4.2} it follows that $q(x) \equiv 0$.

\textbf{Remark 3.} We will get the classical Ambarzumyan's theorem, if we take $\alpha = \pi / 2$.

\

\section{Appendix. Analogues of Theorem \ref{thm1.2}}\label{sec5}

Consider $L(q, \alpha, \beta)$ problem.
Let $\psi(x,\mu)$ be a solution of the equation \eqref{eq1.1}, which satisfies the initial conditions 
\begin{equation}\label{eq5.1}
\psi(\pi,\mu)=1, \quad  \psi'(\pi,\mu)=-\cot\beta.
\end{equation}

The eigenvalues $\mu_n=\mu_n(q, \alpha, \beta)$, $n \geq 0$, are the solutions of the equation
\[
\Phi(\mu) := \varphi(\pi, \mu)\cot \beta+\varphi'(\pi, \mu)=0,
\]
or of the equation
\[
\Psi(\mu) := \psi(0, \mu)\cot\alpha+\psi'(0, \mu)=0.
\]
$\Phi(\mu)$ and $\Psi(\mu)$ are called characteristic functions for the operator $L(q, \alpha, \beta)$.
In \cite{Harutyunyan2} it is proved, that characteristic functions and their derivatives are uniquely determined only from their zeros, i.e. from eigenvalues $\{\mu_n\}_{n=0}^{\infty}$.
It is easy to see that the functions $\psi(x, \mu_n)$, $n \geq 0$, are the eigenfunctions,
corresponding to the eigenvalue $\mu_n$.
The squares of the $L^2$-norm of these eigenfunctions:
\begin{equation*}
b_n = b_n(q,\alpha,\beta) :=\int_0^{\pi} |\psi(x, \mu_n)|^2 dx, \quad   n \geq 0,
\end{equation*}
are called norming constants.

Since all the eigenvalues of $L(q, \alpha, \beta)$ are simple, then there exist constants $\kappa_n=\kappa_n(q, \alpha, \beta)$, $n \geq 0$, such that
\begin{equation}\label{eq5.2}
\varphi(x, \mu_n) = \kappa_n \psi(x, \mu_n).
\end{equation}
The theorem of uniqueness of Harutyunyan (see \cite{Harutyunyan}) states:
\begin{theorem}\label{thm5.1}
	If
	\begin{gather*}
	\mu_n(q, \alpha, \beta) = \mu_n(q_0, \alpha_0, \beta_0), \\
	\kappa_n(q, \alpha, \beta) = \kappa_n(q_0, \alpha_0, \beta_0),
	\end{gather*}
	for all $n \geq 0$, then $\alpha = \alpha_0$, $\beta = \beta_0$ and $q(x) = q_0(x)$ almost everywhere.
\end{theorem}

From \eqref{eq1.4}, \eqref{eq5.1} and \eqref{eq5.2} it follows
\begin{equation}\label{eq5.3}
\kappa_n = \varphi(\pi, \mu_n) = \psi^{-1}(0, \mu_n).
\end{equation}

There is a relationship between norming constants and characteristic functions (see, e.g. \cite{JodLev}, \cite{Harutyunyan}):
\begin{equation}\label{eq5.4}
a_n = |\varphi(\pi, \mu_n)| |\dot{\Phi}(\mu_n)| , 
\end{equation}
\begin{equation}\label{eq5.5}
b_n = |\psi(0, \mu_n)| |\dot{\Psi}(\mu_n)| , 
\end{equation}
where the dot over $\Phi$ (or over $\Psi$) denotes the derivative of $\Phi(\mu)$ with respect to $\mu$.
From equations \eqref{eq5.3} and \eqref{eq5.4} we obtain:
\begin{equation}\label{eq5.6}
a_n = |\kappa_n||\dot{\Phi}(\mu_n)|
\end{equation}

Consider two isospectral operators $L(q, \alpha, \beta)$ and $L(q_0, \alpha_0, \beta_0)$. 
Formulae, analogues to \eqref{eq2.4} and \eqref{eq2.5}, can be obtained for $\kappa_n$:
\begin{equation}\label{eq5.7}
\cot \alpha = \cot \alpha_0 + \sum_{n=0}^{\infty} \cfrac{1}{|\dot\Phi (\mu_n^0)|} \left( \cfrac{1}{|\kappa_n|} - \cfrac{1}{|\kappa_n^0|} \right),
\end{equation}
\begin{equation}\label{eq5.8}
\cot \beta = \cot \beta_0 + \sum_{n=0}^{\infty} \cfrac{|\kappa_n^0| - |\kappa_n|}{|\dot\Phi (\mu_n^0)|}.
\end{equation}

From Theorem \ref{thm5.1} and formulae \eqref{eq2.3}, \eqref{eq5.7}, \eqref{eq5.8}, new statement, similar to Theorem \ref{thm1.2}, can be proven for $\kappa_n$:

\begin{theorem}\label{thm5.2}
	Let $q' \in L^2_\mathbb{R}(0, \pi)$.
	If
	\begin{gather*}
	\mu_n(q, \alpha_0, \beta) = \mu_n(q_0, \alpha_0, \beta_0),\\
	|\kappa_n(q, \alpha_0, \beta)| \geq |\kappa_n(q_0, \alpha_0, \beta_0)|,
	\end{gather*}
	for all $n \geq 0$, then $\beta = \beta_0$ and $q(x) \equiv q_0(x)$.
\end{theorem}

\textbf{Remark 4.} Instead of $\alpha = \alpha_0$ we can fix $\beta = \beta_0$ and/or replace the inequality sign ($"\geq"$) with less then or equal sign ($"\leq"$).
Even so, the result is valid.
Similar theorems can be proven for $\varphi(\pi, \mu_n)$.

\textbf{Remark 5.} Since the uniqueness theorem of Marchenko is also true for norming constants $b_n$, taking into consideration the relations \eqref{eq5.2}, \eqref{eq5.3} and \eqref{eq5.5}, analogues to Theorem \ref{thm1.2} can be proven for $\psi(0, \mu_n)$ and $b_n$.

\

\

\textbf{\Large{Declarations}}

\

\textbf{Acknowledgments} 

The author would like to thank the referees for their helpful comments and suggestions.
The author is also grateful to professor T.N. Harutyunyan for valuable remarks and discussions.

\

\textbf{Funding} 

This work was supported by the RA MES State Committee of  Science, in the frames of the research project No.15T-1A392.

\

\textbf{Competing interests} 

The author declares that he has no competing interests.

\

\textbf{Author's contributions}

The author read and approved the final manuscript.

\

\textbf{Author's information}

Yerevan State University, Alex Manoogian 1, 0025, Yerevan, Armenia.

\

\end{document}